\documentclass[12pt]{amsart}
\usepackage{amsmath}
\usepackage{amssymb}
\usepackage{amsthm}
\usepackage{mathrsfs}
\usepackage{fullpage}
\usepackage{color}
\usepackage[all]{xy}
\usepackage[colorlinks,dvips]{hyperref}

\UseTips

\newcommand{\agn}{\mathscr{A}_{g,N}}

\newcommand{\codim}{\operatorname{codim}}

\newcommand{\er}{\mathbb{E}_\rho}
\newcommand{\erb}{\bar{\mathbb{E}}_\rho}

\newcommand{\fpbar}{\overline{\mathbb{F}}_p}
\newcommand{\fptwo}{\mathbb{F}_{p^2}}
\newcommand{\GL}{\operatorname{GL}}

\newcommand{\GU}{\operatorname{GU}}
\renewcommand{\H}{\operatorname{H}}
\newcommand{\Hom}{\operatorname{Hom}}
\newcommand{\id}{\mathscr{I}_{\Delta}}
\newcommand{\is}{\mathscr{I}_{\Sigma}}
\newcommand{\into}{\hookrightarrow}
\newcommand{\longto}{\longrightarrow}
\newcommand{\kr}{\mathscr{K}_\rho}
\newcommand{\oh}{\mathscr{O}}
\newcommand{\om}{\underline{\omega}}
\newcommand{\R}{\operatorname{R}}
\newcommand{\Res}{\operatorname{Res}}
\newcommand{\shom}{\mathscr{H}\!om}
\newcommand{\Spec}{\operatorname{Spec}}
\newcommand{\ssk}{\mathscr{S}\!\mathscr{S}_{\! k}}
\newcommand{\ssw}{\mathscr{S}\!\mathscr{S}_{\! w}}
\newcommand{\sr}{\mathscr{S}_{\!\rho}}
\newcommand{\zed}{\mathbb{Z}}

\newtheorem{thm}{Theorem}
\newtheorem{cor}[thm]{Corollary}
\newtheorem{prop}[thm]{Proposition}
\newtheorem{fact}[thm]{Fact}
\theoremstyle{definition}
\newtheorem*{ack}{Acknowledgements}
\newtheorem{ex}{Example}
\theoremstyle{remark}
\newtheorem*{rem}{Remark}

\begin{document}
\title{All Siegel Hecke eigensystems (mod $p$) are cuspidal}
\author[A. Ghitza]{Alexandru Ghitza}
%\date{\today}
\address{
  Department of Mathematics and Statistics\\
  McGill University\\
  805 Sherbrooke St. West\\
  Montr\'eal, Qu\'ebec\\
  H3A 2K6, CANADA
}
\email{aghitza@alum.mit.edu}

\begin{abstract}
  Fix integers $g\geq 1$ and $N\geq 3$, and a prime $p$ not dividing $N$.
  We show that the systems of Hecke eigenvalues occurring in the spaces
  of Siegel modular forms (mod $p$) 
  of dimension $g$, level $N$, and varying weight, are the same as the
  systems occurring in the spaces of Siegel \emph{cusp forms} with the
  same parameters and varying weight.  In particular, in the case
  $g=1$, this says that the Hecke eigensystems (mod $p$) coming from
  classical modular forms are the same as those coming from cusp
  forms.  The proof uses both the main theorem of~\cite{Ghitza2004a}
  and a modification of the techniques used there, namely restriction
  to the superspecial locus.
\end{abstract}

\maketitle

\section{Introduction}
This paper is concerned with the systems of Hecke eigenvalues coming
from modular forms in positive characteristic.  The main result is
that imposing the condition of cuspidality has no effect on the set of
eigensystems that can be obtained, at least if we allow ourselves to
change the weight of the form.  We prove this in the context of Siegel
modular forms, but the method should apply to forms coming from other 
Shimura varieties $Y$ of PEL type, i.e. those arising as moduli spaces
of abelian varieties with specified polarizations, endomorphisms, and
level structures.  For the reader's convenience, here
are the key properties of $Y$ that we use: $Y$ should have an
arithmetic Satake compactification $X$, normal and of finite type over
$\Spec\zed[\frac{1}{N}]$ for an appropriate $N$, containing $Y$ as a
dense open subscheme, and such that the superspecial locus of
$Y\otimes\fpbar$ is non-empty.  Moreover, the Hodge
line bundle $\om$ on $Y$ should extend to an ample line bundle on
$X$.  According to \S{}V.0 of~\cite{Faltings1990}, most of these
properties
(except perhaps for the one regarding the non-emptiness of the
superspecial locus) hold for Shimura varieties of PEL
type.

In \S{}\ref{sect:elliptic} we prove the main result in the case of
elliptic modular forms ($g=1$).  We have two reasons for doing this:
first, it will give an idea of the proof of the general case unclouded
by technical complications; second, the elliptic case is simple enough
as to allow us to give an effective version of our result, something
that we cannot accomplish in general.

The rest of the paper deals with the case $g>1$.  In
\S{}\ref{sect:satake} we review properties of the arithmetic Satake
compactification and use them to give our definition of Siegel cusp
forms.  In the process we give a proof of the K\"ocher principle for
the arithmetic Satake compactification, which basically says that any
Siegel modular form extends to the Satake boundary.  This is
undoubtably known to the experts, but it does not seem to have ever
been written down in this setting.  Our proof of the K\"ocher
principle should apply to any Shimura variety of PEL type whose
boundary inside the Satake compactification has everywhere codimension
at least $2$.

In \S{}\ref{sect:superspecial} we
review the definition and properties of superspecial forms.
Section~\ref{sect:proof} contains the proof of our main result.
Finally, in \S{}\ref{sect:chai-faltings} we show that our notion of
Siegel cusp form (defined using the Satake 
compactification) agrees with that introduced by Chai-Faltings
in~\cite{Faltings1990} (based on the toroidal compactifications).

\begin{ack}
  This paper originated in a conversation between the author and
  A.~J.~de~Jong, whom we thank heartily.  Many thanks go to
  A.~Archibald for pointing out a serious bug in an earlier version
  of the paper, to B.~Conrad and P.~Kassaei for helping with the proof
  of Theorem~\ref{thm:extend}, and to E.~Goren and F.~Oort for useful
  clarifications and suggestions.  It is also a pleasure to thank the
  referee for the detailed comments and corrections.  This research
  was funded by the Centre Interuniversitaire en Calcul Math\'ematique
  Alg\'ebrique (CICMA) and the Fonds Qu\'ebecois de la Recherche sur
  la Nature et les Technologies (FQRNT).
\end{ack}

\section{The elliptic case: $g=1$}\label{sect:elliptic}
Let $p$ be a prime, and let $N\geq 3$ be an integer not divisible by
$p$.  We denote by $M_k(N)$ the space of modular forms (mod $p$) of
weight $k$ and level $\Gamma(N)$, and by $S_k(N)$ the subspace of cusp
forms. 
We are interested in the action of Hecke operators $T_\ell$
($\ell\nmid pN$) on the spaces $M_k(N)$ and $S_k(N)$.  There is a
Hecke-equivariant injection
\begin{equation*}
  M_k(N)\into M_{k+p-1}(N)
\end{equation*}
given by multiplication by the Hasse invariant $A$.  The
\emph{filtration} of $f$ is the smallest
integer $w$ such that there exists $\tilde{f}\in M_w(N)$ with
\begin{equation*}
  f=\tilde{f}\cdot(\text{some nonnegative power of }A).
\end{equation*}

In this section we apply Serre's ideas from~\cite{Serre1996} to obtain
a proof of the following
\begin{thm}\label{thm:elliptic}
  Fix a prime $p$ and a level $N\geq 3$, $p\nmid N$.  Then all Hecke
  eigensystems (mod $p$) are cuspidal.  More precisely, let $\Phi$ be
  the eigensystem associated to some $f\in M_k(N)$, and let $w$ be the
  filtration of $f$.  
  \begin{enumerate}
  \item Suppose $p>2$, or $p=2$ and $w>0$.  Then there exists
    $f^\prime\in S_w(N)$ or 
    $S_{w+p^2-1}(N)$ such that $\Phi$ is associated to $f^\prime$.
    Moreover, the first situation (existence of such $f^\prime\in
    S_w(N)$) occurs if and only if $f\in S_k(N)$. 
  \item If $p=2$ and $w=0$, then there exists $f^\prime\in S_6(N)$
    such that $\Phi$ is associated to $f^\prime$.
  \end{enumerate}
\end{thm}
\begin{rem}
  Let $N=3$, $p=2$, $w=0$.  Then $S_0(3)=S_3(3)=0$, so in the
  situation of (b) we cannot do as well as in (a).
\end{rem}
\begin{proof}
  We fix the level $N$ and often drop it from our notation.

  Let $\om$ be the Hodge (line) bundle on $X=X(N)\otimes\fpbar$, so
  that $M_k=\H^0(X,\om^{\otimes k})$, and let $A\in M_{p-1}$
  be the Hasse invariant.  Let $\ssk$ denote the cokernel of
  multiplication by $A$, i.e.
  \begin{equation*}
    0\longto\om^{\otimes (k-p+1)}\xrightarrow{\times A}
    \om^{\otimes k}\longto\ssk\longto 0.
  \end{equation*}
  Define
  \begin{equation*}
    SS_k:=\H^0(X,\ssk)=\H^0(\Sigma,\ssk|_\Sigma),
  \end{equation*}
  where $\Sigma$ denotes the supersingular locus of $X$.  Using the
  fact that every supersingular elliptic curve over $\fpbar$ has a
  canonical $\fptwo$-structure, one easily sees that
  $SS_k=SS_{k+p^2-1}$ for all $k$.  Also $SS_k$ has a natural (away
  from $Np$) Hecke action, coming from the fact that $\ell$-isogenies
  preserve supersingularity (if $\ell\nmid Np$).

  Serre shows (in particular) that if the eigensystem $\Phi$ is
  associated to $f\in M_k$ and $w$ is the filtration of $f$, then
  $\Phi$ also occurs in $SS_w$.  So it is enough for us to show
  that any eigensystem $\Phi$ occurring in $SS_w$ also occurs in
  $S_w$ or $S_{w+p^2-1}$.  To see this, let $\Delta$ denote the divisor of
  cusps on $X$ and consider the exact sequence
  \begin{equation*}
    0\longto\om^{\otimes (w-p+1)}(-\Delta)\xrightarrow{\times A}
    \om^{\otimes w}(-\Delta)\longto\ssw(-\Delta)\longto 0.
  \end{equation*}
  Taking global sections, we get
  \begin{equation}\label{eq:long_ex}
    0\longto S_{w-p+1}\xrightarrow{\times A} S_w\longto SS_w\longto
    \H^1(X,\om^{\otimes (w-p+1)}(-\Delta)),
  \end{equation}
  where we used the identifications
  \begin{equation*}
    \H^0(X,\ssw(-\Delta))=\H^0(\Sigma,\ssw(-\Delta)|_\Sigma)=
    \H^0(\Sigma,\ssw|_\Sigma)=SS_w
  \end{equation*}
  since $\Delta\cap\Sigma=\emptyset$.

  The canonical sheaf of $X$ is $\Omega^1_X$, which is isomorphic
  (by Kodaira-Spencer) to $\om^{\otimes 2}(-\Delta)$.  So by Serre
  duality we have
  \begin{equation*}
    \H^1(X,\om^{\otimes m}(-\Delta))\cong
    \H^0\left(X,\left(\om^{\otimes
          m}(-\Delta)\right)^\vee\otimes\om^{\otimes
        2}(-\Delta)\right)^\vee
    \cong M_{2-m}^\vee,
  \end{equation*}
  where $\cdot^\vee$ denotes the dual.  So if $m>2$ we conclude that
  $\H^1(X,\om^{\otimes m}(-\Delta))=0$.  In particular, if $w>p+1$ we
  know that the map $S_w\to SS_w$ from (\ref{eq:long_ex}) is
  surjective, and therefore $\Phi$ occurs in $S_w$.  

  If $p=2$ and $w=0$, we have $w+2(p^2-1)=6>3=p+1$, so the map $S_6\to
  SS_6=SS_0$ is surjective, and therefore $\Phi$ occurs in $S_6$.
  This settles (b).

  In the situation of (a), $w+p^2-1>p+1$ so the map 
  \begin{equation*}
    S_{w+p^2-1}\longto SS_{w+p^2-1}=SS_w
  \end{equation*}
  is surjective, and therefore $\Phi$ occurs in
  $S_{w+p^2-1}$. 

  It remains to prove the last statement of (a).  Saying that $w$ is the
  filtration of $f$ implies that if we put $n=(k-w)/(p-1)$ then there
  exists $\tilde{f}\in M_w$ such that $A^n\tilde{f}=f$.  Since the
  divisor $\Sigma$ of $A$ is disjoint from the cusps $\Delta$, we have
  that $f\in S_k$ if and only if $\tilde{f}\in S_w$.  So on one hand
  if $f\in S_k$ then $\tilde{f}\in S_w$ and $f$ and $\tilde{f}$ have
  the same eigensystem, so we may put $f^\prime=\tilde{f}$ and we are
  done.  Conversely, suppose $f^\prime\in S_w$, then $f^\prime$ and
  $\tilde{f}$ have the same eigensystem, hence also the same Fourier
  coefficients.  Moreover, they have the same weight $w$, so by the
  $q$-expansion principle $\tilde{f}=f^\prime$.  So $\tilde{f}\in
  S_w$, therefore $f\in S_k$.
\end{proof}

%\begin{rem}
%  We just proved that any non-cuspidal eigenform $f$ of
%  filtration $w$ has the same eigensystem as a cuspidal eigenform of
%  weight $w+p^2-1$.  Here is a shorter proof of this statement: as
%  before, let $n=(k-w)/(p-1)$ and let $\tilde{f}\in S_w$ be such that
%  $A^n\tilde{f}=f$.  Let $f^{\prime}=\theta^{p-1}(\tilde{f})$, where
%  $\theta$ is the operator $M_k\to S_{k+p+1}$ whose commutation
%  relation with the Hecke operators is $T_\ell\circ\theta=\ell\theta\circ
%  T_\ell$.  Then $f^{\prime}\in S_{w+p^2-1}$ has the same eigensystem
%  as $f$.  Note that this argument also works in the pesky special case
%  $w=0$, $p=2$, where the above proof fails.
%\end{rem}

\begin{cor}
  If $f\in M_k(N)$ is an eigenform and has filtration $w>p+1$, then
  $f$ is a cusp form. 
\end{cor}
\begin{proof}%[First proof]
  This follows from the proof of Theorem~\ref{thm:elliptic}: if
  $w>p+1$ then the restriction-to-$\Sigma$ map $S_w\to SS_w$ is
  surjective, so there exists $f^\prime\in S_w$ with the same
  eigensystem as $f$.  But then the last statement of the Theorem
  tells us that $f$ is a cusp form.
\end{proof}
%\begin{proof}[Second proof]
%  This result also follows easily from a well-known fact, namely that
%  every eigensystem is a twist of an eigensystem in weight at most
%  $p+1$.  This means that if $\Phi$ is the eigensystem of $f$, then
%  there exists $g\in M_{k^\prime}$ such that $f=\theta^a(g)$ for some
%  integer $a\geq 0$.  [I want to say that $f$ is then a cusp form, but
%  this is only true if $a>0$.]
%\end{proof}

We conclude this section with some explicit numerical examples, in
which the cusp eigenforms we exhibit are taken from W.~Stein's
database~\cite{MFD}.
\begin{ex}
  Let $p=5$, $N=1$, and let $f=E_4$ be the Eisenstein series of weight
  $4$.  We know that $f$ is a Hecke eigenform with eigensystem
  \begin{equation*}
    (1+\ell^3)_{\ell\neq 5}=(4,3,4,2,3,4,0,3,0,2,4,\dots)
  \end{equation*}
  Of course, $f$ is nothing but the Hasse invariant (mod $5$), so its
  filtration is $w=0$.  Since $f$ is not a cusp form,
  Theorem~\ref{thm:elliptic} predicts the existence of a cuspidal
  eigenform $f^\prime$ of weight $w+p^2-1=24$ with the same
  eigensystem.  Indeed, there is a cusp eigenform (mod $5$) of weight
  $24$ with $q$-expansion 
  {\small
    \begin{eqnarray*}
      f^\prime(q)&=&q + 4q^2 + 3q^3 + 3q^4 + 2q^6 + 4q^7 + 2q^9 +
      2q^{11} + 4q^{12} + 3q^{13} + q^{14} + q^{16} + 4q^{17} +
      3q^{18}\\ 
      &&+ 2q^{21} + 3q^{22} + 3q^{23} + 2q^{26} + 2q^{28} + 2q^{31} +
      4q^{32} + q^{33} + q^{34} + q^{36} + 4q^{37} + O(q^{38}).
    \end{eqnarray*}
  }
% This is just $\Delta(2E_4^3+4E_6^2)$ by comparing $q$-expansions.
  Similarly, $E_6$ has filtration $6$ and eigensystem
  \begin{equation*}
    (1+\ell^5)_{\ell\neq 5}=(3,4,3,2,4,3,0,4,0,2,3,\dots).
  \end{equation*}
  There exists a cusp eigenform (mod $5$) of weight $30$ with
  $q$-expansion
    {\small
      \begin{eqnarray*}
        &&q + 3q^2 + 4q^3 + 2q^4 + 2q^6 + 3q^7 + 3q^9 + 2q^{11} +
        3q^{12} + 4q^{13} + 4q^{14} + q^{16} + 3q^{17} + 4q^{18}\\
        &&+ 2q^{21} +
        q^{22} + 4q^{23} + 2q^{26} + q^{28} + 2q^{31} + 3q^{32} + 3q^{33} +
        4q^{34} + q^{36} + 3q^{37} + O(q^{38}).
      \end{eqnarray*}
    }
% This is $\Delta(3E_4^3E_6+3E_6^3)$.
\end{ex}
\begin{ex}
  Let $p=7$, $N=1$.  We consider the Eisenstein series of weights $4$,
  $6$, and $8$:
  \begin{enumerate}
  \item $f_1=E_4$ has filtration $4$ and eigensystem
    \begin{equation*}
      (1+\ell^3)_{\ell\neq 7}=(2,0,0,2,0,0,0,2,2,0,2,\dots)
    \end{equation*}
    There is a cusp eigenform (mod $7$) of weight $52$ with
    $q$-expansion
    {\small
      \begin{eqnarray*}
        f_1^\prime(q)&=&q + 2q^2 + 3q^4 + 4q^8 + q^9 + 2q^{11} +
        5q^{16} + 2q^{18}+ 4q^{22}\\ 
        &&+ 2q^{23} + q^{25} + 2q^{29} + 6q^{32} + 3q^{36} + 2q^{37} +
        O(q^{38}).
      \end{eqnarray*}
    }
% This is $\Delta(4E_4E_6^6+4E_4^4E_6^4).
  \item $f_2=E_6$ has filtration $0$ and eigensystem
    \begin{equation*}
      (1+\ell^5)_{\ell\neq 7}=(5,6,4,3,0,6,4,5,2,6,5,\dots)
    \end{equation*}
    There is a cusp eigenform (mod $7$) of weight $48$ with
    $q$-expansion
    {\small
      \begin{eqnarray*}
        f_2^\prime(q)&=&q + 5q^2 + 6q^3 + 4q^5 + 2q^6 + q^8 +
        3q^9 + 6q^{10} + 3q^{11} + 3q^{15} + 5q^{16} + 6q^{17} +
        q^{18}+ 4q^{19}\\
        &&+ q^{22} + 5q^{23} + 6q^{24} +
        6q^{25} + 2q^{27} + 2q^{29} + q^{30} + 6q^{31} + 4q^{33} + 2q^{34} + 5q^{37} +
        O(q^{38}) .
      \end{eqnarray*}
    }
  \item $f_3=E_8$ has filtration $8$ and eigensystem
    \begin{equation*}
      (1+\ell^7)_{\ell\neq 7}=(3,4,6,5,0,4,6,3,2,4,3,\dots)
    \end{equation*}
    There is a cusp eigenform (mod $7$) of weight $56$ with
    $q$-expansion
    {\small
      \begin{eqnarray*}
        f_3^\prime(q)&=&q + 3q^2 + 4q^3 + 6q^5 + 5q^6 + q^8 +
        6q^9 + 4q^{10} + 5q^{11} + 3q^{15} + 3q^{16} + 4q^{17} + 4q^{18} +
        6q^{19}\\
        &&+ q^{22} + 3q^{23} + 4q^{24} + 3q^{25} + 5q^{27} + 2q^{29} +
        2q^{30} + 4q^{31} + 6q^{33} + 5q^{34} + 3q^{37} + O(q^{38}). 
      \end{eqnarray*}
    }
  \end{enumerate}
\end{ex}

\section{The arithmetic Satake compactification}\label{sect:satake}
For the remainder of the paper we assume that $g>1$.

Fix an integer $N\geq 3$, and let $\agn$ denote the moduli space of
$g$-dimensional principally polarized abelian varieties with
symplectic level $N$ structure.  

There are several ways to compactify $\agn$; we will work
with the arithmetic Satake (also known as minimal) compactification
$\agn^*$, whose existence and properties are described in Theorem
V.2.5 of~\cite{Faltings1990}.  For now, we just need to know that
$\agn^*$ is a normal scheme, proper and of finite type over
$\Spec\zed[\frac{1}{N}]$, containing $\agn$ as a dense open
subscheme.

Note that in the classical case $g=1$, these
are the usual modular curves $\mathscr{A}_{1,N}=Y(N)$ and
$\mathscr{A}_{1,N}^*=X(N)$.

There is a universal abelian scheme
\begin{equation*}
  \xymatrix{
    {A^{\text{univ}}} \ar^{\pi}[d] & {\Omega^1_{A^{\text{univ}}/\agn}}
    \ar@{~>}[d]\\
    {\agn} & {\mathbb{E}:=\pi_*\Omega^1_{A^{\text{univ}}/\agn}},
  }
\end{equation*}
and we can therefore define the \emph{Hodge bundle} $\mathbb{E}$ on
$\agn$.  It has rank $g$, so given a representation $\rho$ of the
algebraic group $\GL_g$, we can define the \emph{twist} of
$\mathbb{E}$ by $\rho$, by applying $\rho$ to the transition functions
of $\mathbb{E}$.  The result is denoted $\er$ and it has rank equal to
the dimension of $\rho$.

We want to see how much of this can be extended to the minimal
compactification $\agn^*$.  For this we need the following technical
result:

\begin{thm}\label{thm:extend}
  Let $X$ be a locally noetherian scheme which is locally
  of finite type over a quotient $R_0/I$ of a regular ring $R_0$.  Let
  $U$ be an open subset of $X$ such that the complement $Z$ of $U$ has
  everywhere codimension at least $2$ in $X$.  Let $i:U\into X$ denote
  the canonical inclusion.  Let $\mathscr{F}$ be a torsion-free
  coherent $\oh_U$-module.  Then $i_*\mathscr{F}$ is a coherent
  $\oh_X$-module.  If, moreover, $X$ is normal and $\mathscr{F}$ is
  reflexive, so is $i_*\mathscr{F}$, and it is the unique reflexive
  coherent sheaf on $X$ extending $\mathscr{F}$.
\end{thm}
\begin{proof}{}
  We start by noticing that $X$ is \emph{locally embeddable in a
    regular scheme}, i.e. that any point $x\in X$ has an open
  neighborhood isomorphic to a subscheme of a regular scheme.  This
  follows directly from the fact that $X$ is locally of finite type
  over $R_0/I$, i.e. locally embeddable in affine space over $R_0/I$.
  In turn, this is embeddable in a regular scheme, namely affine space
  over the regular ring $R_0$.

  Since $\mathscr{F}$ is torsion-free, the support of $\mathscr{F}$ is
  all of $U$.  Also $U$ is dense in $X$, so $\bar{U}=X$.  Since $Z$ is
  everywhere of codimension at least $2$ in $X$, we have that for any
  irreducible component $U^\prime$ of $\bar{U}=X$, 
  \begin{equation*}
    \codim(Z\cap U^\prime,U^\prime)\geq 2.
  \end{equation*}
  So we may now apply the following result with $n=1$
  to conclude that $i_*\mathscr{F}$ is a coherent $\oh_X$-module:
  \begin{quote}
    {\bf Proposition VIII.3.2 in \cite{SGA2}:}
    Let $X$ be a locally noetherian scheme which is locally embeddable
    in a regular scheme.  Let $U$ be an open subscheme of $X$ and
    let $i:U\to X$ be the canonical embedding.  Let $n\in\zed$, and
    let $\mathscr{F}$ be a coherent Cohen-Macaulay $\oh_U$-module.
    Then the following are equivalent:
    \begin{enumerate}
    \item The sheaf $\R^p i_*\mathscr{F}$ is coherent for all $p<n$.
    \item Let $S$ denote the support of $\mathscr{F}$ and let
      $\bar{S}$ be the closure of $S$ in $X$.  For any irreducible
      component $S^\prime$ of $\bar{S}$, we have
      \begin{equation*}
        \codim(S^\prime\cap(X-U),S^\prime)>n.
      \end{equation*}
    \end{enumerate}
  \end{quote}

  Finally, to prove the statement about reflexivity, we employ an
  argument used by Serre in the complex-analytic category (see
  Proposition 7 of~\cite{Serre1966}).  First notice that
  $i_*\oh_U=\oh_X$: let $V$ be an open subset of $X$ and let
  $f\in\oh_U(U\cap V)$.  We know that $U\cap V$ is dense in $V$ and
  that its complement $Z\cap V$ has everywhere codimension at least $2$ in
  $V$.  So $f$ defines a rational function on $V$; assume that $f$ has
  at least one pole.  Consider the locus $D$ of its poles; by
  assumption $D\subset Z\cap V$, but on the other 
  hand $D$ is a Weil divisor on $V$, so it has codimension one,
  which is absurd since $Z\cap V$ has codimension at least $2$.

  Next we claim that if
  $\mathscr{R}$ is a reflexive sheaf on $X$ then
  \begin{equation*}
    i_*i^*\mathscr{R}=\mathscr{R}.
  \end{equation*}
  To see this, let $\mathscr{D}=\mathscr{R}^\vee$.  Since
  $\mathscr{R}$ is reflexive we have
  $\mathscr{R}=\mathscr{D}^\vee=\shom(\mathscr{D},\oh_X)$, so
  $i^*\mathscr{R}=\shom(i^*\mathscr{D},\oh_U)$.  Therefore
  \begin{equation*}
    i_*i^*\mathscr{R}=i_*\shom(i^*\mathscr{D},\oh_U)=\shom(\mathscr{D},i_*\oh_U)=\shom(\mathscr{D},\oh_X)=\mathscr{D}^\vee=\mathscr{R},
  \end{equation*}
  as claimed.

  Now assume $\mathscr{F}$ is reflexive, and let $\mathscr{G}$ be
  a coherent sheaf on $X$ extending $\mathscr{F}$.  Let
  $\mathscr{G}^{\vee\vee}$ be the bidual of $\mathscr{G}$, then
  $\mathscr{G}^{\vee\vee}$ extends $\mathscr{F}$ and is reflexive.
  Hence
  $\mathscr{G}^{\vee\vee}=i_*i^*\mathscr{G}^{\vee\vee}=i_*\mathscr{F}$, 
  from which we conclude that $i_*\mathscr{F}$ is reflexive and that
  if $\mathscr{G}$ is reflexive then $\mathscr{G}=i_*\mathscr{F}$.
\end{proof}

In particular, we can apply Theorem~\ref{thm:extend} with $X=\agn^*$,
$U=\agn$, $\mathscr{F}=\er$ (we know that the codimension of
$X-U$ in $X$ is $g$, so by our assumptions at least $2$).

We will denote the pushforward $i_*\er$ by $\er^*$.
%According to Exercise II.5.15
%in~\cite{Hartshorne1977}, there exists a \emph{coherent sheaf} $\er^*$
%on $\agn^*$ such that
%\begin{equation*}
%  \er^*|_{\agn}\cong\er.
%\end{equation*}
%Since $\agn$ is a dense open subscheme of $\agn^*$, we conclude that
%this extension $\er^*$ is unique.  
We stress the fact that
$\er^*$ is in general only a coherent sheaf on $\agn^*$, not
necessarily locally free.  This causes some complications in working
with the minimal compactification, but as we shall see they are not
essential.

A notable exception to this caveat, of which we will make crucial use
in our main argument, is the following result (part of Theorem
V.2.5 in~\cite{Faltings1990}): 
\begin{fact}[Chai-Faltings]\label{thm:omega}
  The invertible sheaf $\om:=\mathbb{E}_{\det}$ on $\agn$ extends
  to an invertible sheaf $\om^*$ on $\agn^*$ relatively ample over
  $\Spec\zed$.
\end{fact}

We have the following result:
\begin{prop}[K\"ocher principle]\label{prop:koecher}
  If $g>1$, then for any $\rho$ and for any $\zed[\frac{1}{N}]$-module
  $M$ there is a natural identification 
  \begin{equation*}
    \H^0(\agn^*,\er^*\otimes M)=\H^0(\agn,\er\otimes M).
  \end{equation*}
\end{prop}
\begin{proof}
  This is presumably well-known but we prove it here for lack of a
  reference.  See Theorem~10.14 of~\cite{Baily1966c} for the
  complex-analytic version of a more general result.

  Any $\zed[\frac{1}{N}]$-module $M$ is a direct limit of finitely
  generated $\zed[\frac{1}{N}]$-modules; since cohomology and tensor
  products commute with direct limits, we may safely assume that $M$
  is finitely generated.  Using the classification of finitely
  generated $\zed[\frac{1}{N}]$-modules and the additivity of
  cohomology, we may further reduce to the case where
  $M=\zed[\frac{1}{N}]$ or $M=\zed/n\zed$ for some integer $n$ coprime
  to $N$.  In both cases $M$ is actually a ring, which we denote by
  $R$.

  Let $X=\agn^*\otimes R$ and $Y=\agn\otimes R$.  Since
  $\er^*=i_*\er$, we have
  \begin{equation*}
    \H^0(X,i_*\er)=(i_*\er)(X)=\er(Y)=\H^0(Y,\er),
  \end{equation*}
  as desired.
\end{proof}

Given a $\zed[\frac{1}{N}]$-module $M$, the space of \emph{Siegel
  modular forms} of weight $\rho$ and level $N$ with coefficients in
$M$ is  
\begin{equation*}
  M_\rho(M):=\H^0(\agn,\er\otimes M).
\end{equation*}
(We do not include $N$ in the notation since we will always work
with a fixed $N$.)

%\subsection{Siegel cusp forms}\label{sect:cusp}
We define the \emph{cusps} to be the boundary
\begin{equation*}
  \Delta:=\agn^*-\agn.
\end{equation*}
We want to define a notion of cusp form in this setting.  We start
with the short exact sequence of $\oh_{\agn^*}$-modules that defines
the ideal sheaf $\id$ of $\iota:\Delta\into\agn^*$:
\begin{equation*}
  0\longto\id\longto\oh_{\agn^*}\longto\iota_*\oh_\Delta\longto 0,
\end{equation*}
and we tensor it with $\er^*$; since $\er^*$ is not necessarily
locally free, we only get
\begin{equation*}
  \id\otimes\er^*\longto\er^*\longto\er^*|_\Delta\longto 0.
\end{equation*}
Define
\begin{equation*}
  \sr:=\ker\left(\er^*\longto\er^*|_\Delta\right).
\end{equation*}
In other words, for any open $U\subset\agn^*$, $\sr(U)$ consists of
the sections of $\er$ over $U$ that vanish at the cusps $\Delta$.

It is then natural to define the space of \emph{Siegel cusp forms} of
weight $\rho$ and level $N$ with coefficients in a
$\zed[\frac{1}{N}]$-module $M$ to be
\begin{equation*}
  S_\rho(M):=\H^0(\agn^*,\sr\otimes M).
\end{equation*}
Note that as a result of this definition and of the K\"ocher
principle, $S_\rho(M)$ is a subset of $M_\rho(M)$.

\section{Hecke eigensystems (mod $p$)}
We now fix a prime $p$ not dividing $N$, and set
\begin{equation*}
  U:=\agn\otimes\fpbar,\quad X:=\agn^*\otimes\fpbar,\quad
  M_\rho:=M_\rho(\fpbar),\quad S_\rho:=S_\rho(\fpbar).
\end{equation*}
There is a Hecke action on $M_\rho$ given by the Hecke operators
corresponding to the primes not dividing $Np$.  They are essentially
induced by isogenies of degree coprime to $Np$ (for their exact
definition, see \S{}2.2.2 and \S{}3.2.1 of~\cite{Ghitza2004a}).  We
denote by $\mathscr{H}$ the $\zed$-algebra generated by these
operators.  It is known to be commutative (see Satz~IV.1.13
of~\cite{Freitag1983}). 

If $V$ is any $\fpbar$-vector space with an action of $\mathscr{H}$,
an element $v\in V$ which is a common eigenvector for $\mathscr{H}$
defines an algebra homomorphism $\Phi:\mathscr{H}\to\fpbar$ given by
\begin{equation*}
  Tv=\Phi(T)v,\quad\text{for all }T\in\mathscr{H}.
\end{equation*}
This $\Phi$ is called the \emph{eigensystem} associated to $v$.

With this terminology, a \emph{Hecke eigensystem} (mod $p$) is one
associated to an element of $M_\rho$.  The action of $\mathscr{H}$
restricts to $S_\rho$; this follows from the fact that $S_\rho$ is the
space of global sections of the coherent sheaf $\mathscr{S}_\rho$ on
$X$ (see \S{}2.2.2 of~\cite{Ghitza2004a}).  We say that $\Phi$ is
\emph{cuspidal} if it is associated to an element of $S_\rho$.

We can now state our main result:
\begin{thm}\label{thm:main}
  Fix the characteristic $p>0$, the dimension $g\geq 2$, and the level
  $N\geq 3$, $p\nmid N$. 
  Then all Hecke eigensystems (mod $p$) are cuspidal.  That is, for
  any $\Phi$ associated to some $f\in M_\rho$ there exists some $f'\in
  S_{\rho'}$ such that $\Phi$ is associated to $f'$.
\end{thm}

The idea of the proof is this: we use a result of~\cite{Ghitza2004a}
saying that any eigensystem (mod $p$) is superspecial (see the next
section for the definition), and then we show that any superspecial
eigensystem is cuspidal.

\section{Superspecial forms}\label{sect:superspecial}
A $g$-dimensional abelian variety $A$ over $\fpbar$ is said to be
\emph{superspecial} if
\begin{equation*}
  \dim_{\fpbar}\Hom(\alpha_p,A)=g,
\end{equation*}
where $\alpha_p$ is the kernel of multiplication by $p$ on the
additive group $\mathbb{G}_a$ over $\fpbar$.  Equivalently, $A$ is
$\fpbar$-isomorphic to $E^g$, where $E$ is any supersingular
elliptic curve.

Let $\Sigma\subset X$ denote the set of superspecial points.  It has
a number of remarkable properties, including
\begin{itemize}
\item It is finite.
\item It is closed under isogenies of degree coprime to $Np$.
\item Any superspecial $A$ has a canonical and functorial
  $\fptwo$-structure (see Proposition 6 in~\cite{Ghitza2004a}); in
  particular it 
  makes sense to talk about the space of $\fptwo$-rational
  differentials on $A$, and it turns out that a principal polarization
  on $A$ induces a natural hermitian form on this space.  Thus if we
  are interested in hermitian bases, the change-of-basis group is  
  \begin{equation*}
    \GU_g(\fptwo):=\{M\in\GL_g(\fptwo):{}^t\overline{M}M=
    \gamma(M)I\text{ for some }\gamma(M)\in\fptwo^\times\},
  \end{equation*}
  where the ``conjugation'' $\overline{\cdot}:\fptwo\to\fptwo$ is
  $a\mapsto \overline{a}=a^p$.
\end{itemize}
This suggests the following definition.  Given a finite-dimensional
$\fpbar$-representation 
\begin{equation*}
  \tau:\GU_g(\fptwo)\longto W,
\end{equation*}
we set
\begin{multline*}
  SS_\tau:=\{f:[A,\lambda,\alpha,\eta]\longto W\text{ such that}\\
  f([A,\lambda,\alpha,M\eta])=\tau(M)^{-1}f([A,\lambda,\alpha,\eta])
  \text{ for all }M\in\GU_g(\fptwo)\},
\end{multline*}
where $[A,\lambda,\alpha,\eta]$ denotes the $\fpbar$-isomorphism class
of the quadruple, and
\begin{itemize}
\item $(A,\lambda)$ is a superspecial principally polarized abelian
  variety over $\fpbar$;
\item $\alpha$ is a symplectic level $N$ structure on $(A,\lambda)$;
\item $\eta$ is a hermitian basis of invariant $\fptwo$-rational
  differentials on $(A,\lambda)$.
\end{itemize}

We refer to $SS_\tau$ as the space of \emph{superspecial forms} of
weight $\tau$.  The aforementioned properties of $\Sigma$ imply that
$SS_\tau$ admits an action of the Hecke algebra $\mathscr{H}$ and that
it has the following periodicity property:
\begin{equation*}
  SS_{\tau\otimes\det^{\otimes p^2-1}}=SS_\tau\quad\text{for all }\tau.
\end{equation*}
An eigensystem $\Phi$ associated to some $f\in SS_\tau$ is said to be
\emph{superspecial}.

\section{Proof of the main result}\label{sect:proof}
We now prove Theorem~\ref{thm:main}.  

It is part of the proof of Theorem~28 in~\cite{Ghitza2004a}
(more precisely, the first paragraph on page 380 loc. cit.) that any
eigensystem (mod $p$) is superspecial.  Therefore it suffices to show
that any superspecial eigensystem is cuspidal.

Let $\is$ be the ideal sheaf of $j:\Sigma\into X$; it is defined by
the short exact sequence of $\oh_X$-modules
\begin{equation*}
  0\longto \is\longto \oh_X\longto j_*\oh_\Sigma\longto 0.
\end{equation*}
Upon tensoring with the coherent sheaf $\sr$ introduced
towards the end of \S{}\ref{sect:satake}, we get
\begin{equation}\label{eq:srs}
  \is\otimes\sr\longto\sr\longto\sr|_\Sigma\longto 0.
\end{equation}

We can easily pass from $\sr|_\Sigma$ to $\er|_\Sigma$; since
restriction to $U$ is exact and $\Delta\cap U=\emptyset$, the short
exact sequence defining $\sr$
\begin{equation*}
  0\longto\sr\longto\er^*\longto\er^*|_\Delta\longto 0
\end{equation*}
gives an isomorphism $\sr|_U\cong\er^*|_U$.  In particular,
$\sr|_\Sigma\cong\er^*|_\Sigma\cong\er|_\Sigma$, the latter
isomorphism coming from $\Sigma\cap\Delta=\emptyset$.

Therefore the surjective map from the sequence~(\ref{eq:srs}) gives
\begin{equation*}
  \sr\longto\er|_\Sigma\longto 0.
\end{equation*}
Let $\kr$ denote its kernel:
\begin{equation}\label{eq:kr}
  0\longto\kr\longto\sr\longto\er|_\Sigma\longto 0.
\end{equation}
This yields a long exact sequence
\begin{equation*}
  0\longto\H^0(X,\kr)\longto S_\rho\longto\H^0(\Sigma,\er|_\Sigma)
  \longto\H^1(X,\kr).
\end{equation*}
But it is easily seen that
\begin{equation*}
  \H^0(\Sigma,\er|_\Sigma)=SS_{\Res\rho},
\end{equation*}
where $\Res\rho$ denotes the restriction of $\rho$ to the finite group
$\GU_g(\fptwo)$.

Therefore we have a map (which we think of as restriction of cusp
forms to the superspecial locus)
\begin{equation*}
  r_\rho:S_\rho\to SS_{\Res\rho},
\end{equation*}
which is Hecke-equivariant and whose cokernel is contained in
$\H^1(X,\kr)$.  Recall that Fact~\ref{thm:omega} says that
$\om^*$ is a line bundle; by tensoring the short exact
sequence~(\ref{eq:kr}) with $\om^*$, it is easy to see that
\begin{equation*}
  \mathscr{K}_{\rho\otimes\det}=\kr\otimes\om^*.
\end{equation*}
Since $\om^*$ is an ample line bundle on the projective scheme $X$
over $\fpbar$, we know from a theorem of Serre (Theorem~III.5.2
in~\cite{Hartshorne1977}) that for $k$ large enough we have
\begin{equation*}
  \H^1(X,\mathscr{K}_{\rho\otimes\det^k})=
  \H^1(X,\kr\otimes(\om^*)^{\otimes k})=0.
\end{equation*}
Thus for $k$ large enough we know that $r_{\rho\otimes\det^k}$ is
surjective.

Now suppose we start with a superspecial eigensystem $\Phi$, say
associated to some $f\in SS_\tau$.  By Corollary~27
of~\cite{Ghitza2004a}, we can extend $\tau$ to $\GL_g(\fpbar)$,
i.e. there exists a rational representation
\begin{equation*}
  \rho:\GL_g(\fpbar)\longto\GL(V)\text{ such that
  }\tau\subset\Res\rho.
\end{equation*}
This means that $SS_\tau\subset SS_{\Res\rho}$.  Now by the
periodicity property of $SS_\tau$ we have
\begin{equation*}
  SS_{\Res\rho}=SS_{\Res\rho\otimes\det^{k(p^2-1)}}\text{ for all }k.
\end{equation*}
So we can pick $k$ large enough such that
\begin{equation*}
  r_{\rho\otimes\det^{k(p^2-1)}}:S_{\rho\otimes\det^{k(p^2-1)}}\longto
  SS_{\Res\rho\otimes\det^{k(p^2-1)}}=SS_{\Res\rho}\supset SS_\tau
\end{equation*}
is surjective.  Therefore by a simple linear algebra argument we
conclude that $\Phi$ is associated to some $f'\in
S_{\rho\otimes\det^{k(p^2-1)}}$.

\section{Comparison with cusp forms \`a la
  Chai-Faltings}\label{sect:chai-faltings} 
Siegel cusp forms were already defined in an algebraic-geometric way
by Chai and Faltings (see pp. 144--147 of~\cite{Faltings1990}), at
least in the scalar case, i.e. for $\rho=\det^{\otimes k}$.  Since
their definition is based on the toroidal 
compactifications of $\agn$, it is not immediately obvious that it
agrees with ours, and this last section is devoted to showing this.

Chai and Faltings define so-called arithmetic toroidal compactifications
$\bar{\agn}$ of $\agn$.  These depend on certain combinatorial data,
and have various nice properties summarized in Theorem~IV.6.7
of~\cite{Faltings1990}.  Most importantly, they are moduli spaces and
thus one has a Hodge bundle $\bar{\mathbb{E}}$ and its twisted
versions $\erb$, which we define in the same way as we did
$\mathbb{E}$ and $\er$ in~\S{}\ref{sect:satake}.  In this setting,
Chai and Faltings define Siegel cusp forms with coefficients in a
$\zed[\frac{1}{N}]$-module $M$ to
be
\begin{equation*}
  \H^0(\bar{\agn},\mathcal{I}_{\bar{\Delta}}\otimes\erb\otimes M),\quad
  \text{where }\mathcal{I}_{\bar{\Delta}}\text{ is the ideal sheaf of
  }\bar{\Delta}=\bar{\agn}-\agn.
\end{equation*}
In other words, cusp forms are global sections of $\erb\otimes M$ that
vanish along the boundary $\bar{\Delta}$ of $\bar{\agn}$.  Moreover,
this turns out to be independent of the choice of toroidal
compactification.

The key to comparing the two notions of cusp forms is the following
fact (part of Theorem~V.2.5 of~\cite{Faltings1990}): a toroidal
compactification is related to the minimal compactification by a
canonical morphism $\bar{\pi}:\bar{\agn}\to\agn^*$ restricting to the
identity on the open dense subscheme $\agn$.  Two facts about
$\bar{\pi}$ are important for our purposes:
\begin{itemize}
\item The boundary $\bar{\Delta}$ is the scheme-theoretic preimage of
  $\Delta\subset\agn^*$ under $\bar{\pi}$: this follows from the
  detailed description of the interaction between $\bar{\pi}$ and the
  stratifications of $\bar{\agn}$ and $\agn^*$ (see Theorem~V.2.5(6)
  of~\cite{Faltings1990}).
\item The pullback $\bar{\pi}^*(\er^*)$ is $\erb$: this can be
  seen easily from the fact that both are reflexive coherent sheaves
  on $\bar{\agn}$ extending $\er$ on $\agn$, together with the
  uniqueness argument from the end of Theorem~\ref{thm:extend}.
\end{itemize}

Using these it is easy to obtain the following result, whose proof we
leave to the reader:
\begin{prop}\label{prop:comp}
  The canonical morphism $\bar{\pi}:\bar{\agn}\to\agn^*$ induces via
  pullback an isomorphism 
  \begin{equation*}
    \xymatrix{
      {\bar{\pi}^*:S_\rho(M)} \ar[r]^-{\cong} &
      {\H^0(\bar{\agn},\mathscr{I}_{\bar{\Delta}}\otimes\erb\otimes M)}
    }
  \end{equation*}
  for any $\zed[\frac{1}{N}]$-module $M$.
\end{prop}

%\bibliographystyle{halpha}
%\bibliography{mrl}

\end{document}